\documentclass[11pt]{amsart}
\input psfig.sty

\newcommand{\T}{\nabla}

\newcommand{\rr}{\mathbb{R}}

\theoremstyle{plain}
\newtheorem{thm}{Theorem}[section]
\newtheorem{lem}{Lemma}[section]

\theoremstyle{definition}

\newtheorem{defi}{Definition}[section]

\theoremstyle{remark}
\newtheorem{remark}{Remark}

\newcommand{\1}{\hspace{1mm}}

\newcommand{\5}{\hspace{5mm}}
\newcommand{\6}{\hspace{6mm}}

\newcommand{\8}{\hspace{8mm}}

\title[Poisson's Equation with Interfaces]
{Convergence of Finite Difference Methods for 
Poisson's Equation with Interfaces}
\author{Xu-Dong Liu }
\thanks{Research partially supported by the National Science Foundation:
DMS-9805546 (X.-D.L.) and DMS-9800888 (T.C.S.).} 
\author{Thomas C. Sideris}
\address{Department of Mathematics\\
University of California\\
Santa Barbara, CA 93106}
\email{xdl@math.ucsb.edu, sideris@math.ucsb.edu}
\thanks{AMS subject classification:  65N12, 35J25}

\begin{document}

\date{August 15, 2001}

\begin{abstract}
\noindent
In this paper, a weak formulation of the discontinuous variable
coefficient
Poisson equation with interfacial jumps is studied. The existence, 
uniqueness and regularity of solutions of this problem are obtained. 
It is shown that the application of the Ghost Fluid Method by Fedkiw, 
Kang, and Liu to this problem in \cite{lfk} can be obtained in a 
natural way through 
discretization of the weak formulation.  An abstract framework is given for
proving the convergence of finite difference methods derived from a weak
problem, and as a consequence, the Ghost Fluid Method is proven to be
convergent.

\end{abstract}

\maketitle

\numberwithin{equation}{section}

\section{Introduction}

This paper proves the convergence of the finite difference 
method introduced in \cite{lfk}
for the Poisson equation with discontinuous
coefficients and given interfacial jumps. Based on the Ghost Fluid Method (GFM)
of \cite{ghost}, see also \cite{deflagration}, \cite{XuDong},
the finite difference method is simple, efficient and robust.
One of the novelties, and advantages, of the method is the arm-by-arm splitting 
technique which makes the method in multi-dimension as simple 
as in one-dimension. Another advantage
is that the resulting linear system of the method is the same as the linear system 
obtained from the simplest standard five point stencil finite difference method
for the Poisson equation without discontinuous
coefficients or given interfacial jumps. Therefore the resulting linear system is
symmetric and positive definite and can be efficiently solved.
Inherited from the GFM, this finite difference method captures the
sharp solution profile at the interfaces without smearing.
This is a necessity in the development
of effective schemes for problems involving interfaces.
A variety of other approaches to problems with interfaces
have been taken, \cite{colella}, \cite{LevequeLi}, \cite{fastLi},
\cite{mulder}, \cite{peskin}, \cite{peskin2}.
An important application of the method here is for
Hele-Shaw flow, see \cite{HeleShaw}.

The convergence
proof starts with the formulation 
of the problem, in terms of a uniformly elliptic
bilinear form.
Existence and uniqueness of a weak solution $v$
follow immediately using elementary functional analysis.
The solution space is the same as $H^1_0$ 
except with a different, but equivalent,
inner product induced by the bilinear form.
Discretizing this weak formulation in a 
natural way results a discrete weak problem,
which is equivalent to the finite different method 
in \cite{lfk}. As expected, the
discrete weak formulation inherits the 
structural conditions (uniform boundedness, extension, consistency) in the discrete sense,
hence existence, uniqueness, and uniform boundedness of the
family $v^h$ of discrete solution hold. 

We then provide an abstract framework for the convergence
proof.  In addition to the uniform structural conditions
for the weak problem and its discretizations, we
postulate the existence of a uniformly bounded
family of extension operators $T^h$
from the discrete spaces to the solution space, which
satisfy a strong approximation property.
This natural assumption implies that 
as the grid size $h$ goes to zero the image
of the discrete space fills out the entire
solution space.
Finally, we impose  weak consistency of the
discrete problem with the original problem.
Under these assumptions, the sequence of extended 
weak solutions $T^h(v^h)$ converges weakly to the
solution $v$ in $H^1_0$.
Guided by the abstract framework, we 
establish the converge of the finite difference
scheme of \cite{lfk}.

There is a similarity between the standard convergence proof of finite element 
methods and our approach. Both of them use structural conditions.
For finite element methods, the structural conditions
and weak consistency are inherited directly from the weak
problem for the PDE because the discrete bilinear forms
are obtained by restriction to finite dimensional subspaces.
Cea's Lemma then says that the extended discrete solution
is the closest function in the finite dimensional subspace
to the true solution.  This reduces convergence and
error estimation to a problem in approximation theory.
For finite difference methods, further approximations
are made so that
the discrete problem can not be obtained by restriction to
a finite dimensional subspace.  
In our case, the difference between the finite difference
scheme and the finite element method is that
point values of the coefficients are used rather
than cell averages.  
Because the coefficients in our problem are
discontinuous, the solution lies only in $H^1_0$ and not $H^2\cap H^1_0$,
and therefore, we obtain convergence, but without a rate.

\section{Equations and Weak Formulation}
\setcounter{equation}{0}

Consider a bounded domain, $\Omega\subset\rr^n$, with 
smooth boundary, $ \partial \Omega $.
Let $ \Gamma\subset\Omega$ be a smooth 
interface of co-dimension $n-1$,
represented by the zero level-set of
a smooth function $\phi(x)$,
which is a signed distance function of the interface locally.
We  assume that $\phi$ divides 
$\Omega$ into disjoint subdomains, $ \Omega^- =\{\phi<0\}$ and 
$ \Omega^+=\{\phi>0\}$, with $\partial\Omega^-=\Gamma$. 
Thus, we may write $\Omega=\Omega^+\cup\Omega^-\cup\Gamma$.
The unit normal vector of the interface is 
$ n=\nabla\phi/|\nabla\phi| $, for $\phi(x)=0$, pointing
from $\Omega^-$ to $\Omega^+$.

We seek solutions of the
variable coefficient Poisson equation away from the interface
given by
\begin{subequations}
\begin{equation}
\nabla \cdot (\beta({x})\nabla
u({x}))=f({x}), \quad {x}\in\Omega\setminus\Gamma,
\label{pde(a)}
\end{equation}
in which $ {x} = (x_1,\ldots,x_n) $ denotes the spatial variables and
$ \nabla$
is the gradient operator.  The coefficient $\beta({x})$
is assumed to be a positive definite, symmetric $n\times n$
matrix, the components of which are continuously
differentiable on the closure of each disjoint subdomain, 
$ \Omega^- $ and $ \Omega^+ $, but they may be discontinuous
across the interface $ \Gamma $.  
It follows that there are positive constants $m<M$ with
$m \;Id\leq \beta({x})\leq M\; Id$, where $Id$ stands for
the $n\times n$ identity matrix.
We suppose that on the interface,
$\beta$ assumes the limiting values from within $\Omega^-$.
%The function $\beta({x})$ could be generalized to the 
%case of symmetric matrix.  
The right-hand side $f(x)$ is assumed to lie in $L^2(\Omega)$.
 
Given functions $a$ and $b$ along the interface $ \Gamma $,
we prescribe the jump conditions 
\begin{equation}
\left\{\begin{array}{l}
\left[u\right]_{\Gamma}(x) 
\equiv u^+({x})-u^-({x})
=a({x}),\hspace{8mm}\\
\left[(\beta u)_n\right]_{\Gamma}(x)
\equiv (\beta u)_n^+({x})-(\beta u)_n^-({x})
=b({x}),
       \end{array}\right.
\quad {x}\in\Gamma.
\label{pde(b)}
\end{equation}
Note that $ (\beta u)_n = n \cdot\beta\nabla u$,
and the ``$\pm$'' subscripts refer to limits taken
from within the subdomains $ \Omega^{\pm} $.

Finally, we prescribe boundary conditions
\begin{equation}
\label{pde(c)}
u(x)=g(x), \quad x\in\partial\Omega,
\end{equation}
\end{subequations}
for a given function $g$ on the boundary.

We are going to use the usual Sobolev spaces
$H^1_0(\Omega)$ and $H^1(\Omega)$, but instead 
of the usual inner product we choose 
one which is better suited to our problem:
\begin{equation}
\label{bilinearform}
B[u,v]
=\int\limits_{\Omega}\beta\T u\cdot\T v.
\end{equation}
This induces a norm on $H^1_0(\Omega)$ which is equivalent to the usual
one, thanks to the Poincar\'e inequality and the uniform
bounds for the coefficient matrix.

Let $R_\Gamma$ and $R_{\partial\Omega}$
denote the restriction operators
from $H^1(\Omega)$ to $L^2(\Gamma)$ and $L^2(\partial\Omega)$,
respectively.
Throughout this section,
we shall always assume that our boundary data $a$, $b$
are the restrictions of functions $\widetilde{a}$, 
$\widetilde{b}\in H^1_0(\Omega)$, and that
$g$ is the restriction of a function
$\widetilde{g}\in H^1(\Omega)$:
\begin{equation}
\label{restriction}
a = R_\Gamma(\widetilde{a}),
\quad
b = R_\Gamma(\widetilde{b}),
\quad
\mbox{and}
\quad
g = R_{\partial\Omega}(\widetilde{g}).
\end{equation}
To simplify the notation, from now on we will drop
the tildes.

First let us consider the version of problem of (\ref{pde(a)}),
(\ref{pde(b)}), \eqref{pde(c)}
 with homogeneous boundary conditions for $u$:
\begin{equation}
\left\{\begin{array}{ll}
         \T \cdot \left(\beta({ x})
            \T u({ x})\right)
            =f({ x}), \hspace{8mm}&{ x}\in\Omega\setminus\Gamma\\
         \left[u\right]_{\Gamma}(x)=0,& x\in\Gamma\\
         \left[(\beta u)_n\right]_{\Gamma}(x)=b({ x}),& x\in\Gamma\\
         u( x)=0,
            &{ x}\in\partial\Omega.
      \end{array}\right.
\label{pde_homo_bc}
\end{equation}

\begin{defi}\label{weak_solu_homo_bc}
A function $u \in H_0^1(\Omega)$ is a weak solution
of (\ref{pde_homo_bc}) if it satisfies
\begin{equation}\label{weak_form_homo_bc}
-B[ u,\psi]=
\int\limits_\Omega f\psi
+\int\limits_{\Omega^-}\T\cdot(b\psi{n})\1,
\end{equation}
for all $\psi\in H_0^1(\Omega)$.
\end{defi}

A classical solution of (\ref{pde_homo_bc}),
$u|_{\Omega^\pm}\in {C^2(\overline{\Omega^\pm})}$,
  is easily seen to satisfy
(\ref{weak_form_homo_bc}).
The boundary condition $[u]_\Gamma=0$ is understood with the aid
of  $R_\Gamma$, and it is fulfilled since $u\in H^1_0(\Omega)$.

\begin{thm} \label{exist_uniq_homo_bc}
Given $f\in L^2(\Omega)$ and $b\in H^1_0(\Omega)$, there exists
a unique weak solution of \eqref{weak_form_homo_bc} in $H_0^1(\Omega)$.
\end{thm}

\begin{proof}
The right-hand side of (\ref{weak_form_homo_bc})
\begin{displaymath}
F(\psi)=\int\limits_{\Omega}f\psi
+\int\limits_{\Omega^-}\T\cdot(b\psi{n})
\end{displaymath}
is a continuous linear functional on $H_0^1(\Omega)$.  By the Riesz
representation theorem, there exists a unique $u\in H_0^1(\Omega)$
such that $-B[u,\psi]=F(\psi)$, for all $\psi\in H_0^1(\Omega)$.
\end{proof}

Next we reduce the general case (\ref{pde(a)}), (\ref{pde(b)}),
\eqref{pde(c)}
to the homogeneous case.
We will construct a unique solution of the problem
in the class
\begin{displaymath}
H(a,g)=\{u :
u-g+a\chi_{\Omega^-}\in H^1_0(\Omega)\},
\end{displaymath}
in which $\chi_{\Omega^-}$ is the characteristic function of 
$\Omega^-$.
If $u\in H(a,g)$ then
\[
[u]_\Gamma = a
\quad\mbox{and}\quad
R_{\partial\Omega}(u)=g.
\]
Note that $H_0^1(\Omega)$ can be identified with $H(0,0)$, and so
the following definition is consistent with the previous one.

\begin{subequations}
\begin{defi}\label{weak_solu}
A function $u\in H(a,g)$ is a weak solution
of (\ref{pde(a)}), (\ref{pde(b)}), \eqref{pde(c)}
 if $v=u-g+a\chi_{\Omega^-}$ 
satisfies
\begin{equation}\label{weak_form(a)}
-B[ v,\psi] =F(\psi),
\end{equation}
for all $\psi\in H_0^1(\Omega)$,
where
\begin{align}\label{weak_form(b)}
F(\psi)=&F_1(\psi)+\ldots+F_4(\psi)\\
\nonumber
=&\int\limits_{\Omega}f\psi
+\int\limits_{\Omega}\beta\T g\cdot\T\psi
-\int\limits_{\Omega^-}\beta\T a\cdot\T\psi
+\int\limits_{\Omega^-}\T\cdot(b\psi{n}).
\end{align}
\end{defi}
\end{subequations}

A classical solution of (\ref{pde(a)}), (\ref{pde(b)}),
\eqref{pde(c)}
is necessarily a weak solution.

\begin{thm}\label{exist_uniq}
If $f\in L^2(\Omega)$, $g\in H^1(\Omega)$, and $a$, $b\in H^1_0(\Omega)$,
then there exists a unique weak solution
of (\ref{weak_form(a)}), (\ref{weak_form(b)}), \eqref{pde(c)} in $H(a,g)$.
\end{thm}

\begin{proof}
The right-hand side of (\ref{weak_form(b)})
is a continuous linear functional on $H_0^1(\Omega)$.  By the Riesz
representation theorem, there exists a unique $v\in H_0^1(\Omega)$
such that $-B[ v,\psi]=F(\psi)$, for all $\psi\in H_0^1(\Omega)$.
Define the solution as $u=v+g-a\chi_{\Omega^-}$.
\end{proof}

\begin{remark}
We could replace $\int\limits_{\Omega^-}\T\cdot(b\psi{n})$
by $\int\limits_\Gamma b\psi\,ds$ in (\ref{weak_form(b)}). Then 
the requirement for $b$ is reduced to $b\in L^2(\Gamma)$, which
is a bit weaker than our assumption \eqref{restriction} for $b$.
\end{remark}

\section{Abstract Form of Finite Difference Methods}

\6 In this section an abstract framework is given for proving convergence 
of finite difference methods for the elliptic interface problem 
introduced in the previous section.

For any $h>0$, let $H_0^{1,h}$ be a finite dimensional vector
space with norm $\|\cdot\|_h$.  This space should
be thought of as a discrete approximation of the
Sobolev space $H^1_0(\Omega)$ with grid size measured
by the parameter $h$.  

On each finite dimensional space $H_0^{1,h}$,
we suppose there exists a bounded {\bf extension
operator} $T^h:H_0^{1,h}\to H^1_0(\Omega)$
with the bound
\begin{subequations}
\begin{equation}
\label{bdd_ext}
\|T^h(\psi^h)\|_{H^1_0(\Omega)} \le C_0 \|\psi^h\|_h,
\end{equation}
for all $\psi^h\in H_0^{1,h}$, with $C_0$ independent of $h$.

Furthermore, we assume a {\bf strong approximation property}.
That is, given $\psi\in H^1_0(\Omega)$, there exists a sequence
$\psi^h\in H_0^{1,h}$ such that
\begin{equation}
\label{str_app}
T^h(\psi^h)\to \psi
\quad \mbox{in}\quad H^1_0(\Omega),
\quad \mbox{as}\quad h\to 0.
\end{equation}
\end{subequations}

On each vector space $H_0^{1,h}$, we suppose there is a 
{\bf uniformly bounded family of bilinear forms}
$B^h[\cdot,\cdot]$ such that for every $u^h$, $v^h\in H_0^{1,h}$
\begin{subequations}
\begin{equation}\label{abs_structure(a)}
C_1\|u^h\|^2_{h}\leq B^h[u^h,u^h]\5\mbox{and}\5
\mid B^h[u^h,v^h]\mid\leq C_2\|u^h\|_{h}\|v^h\|_{h},
\end{equation}
for positive constants $C_1$, $C_2$ independent of $h$.
We also assume there exists a {\bf uniformly bounded
family of linear functionals}
$F^h(\cdot)$ on $H_0^{1,h}$ such that for every $\psi^h\in H_0^{1,h}$
\begin{equation}\label{abs_structure(b)}
\mid F^h(\psi^h)\mid\leq C_3\|\psi^h\|_h,
\end{equation}
\end{subequations}
again with a constant $C_3$ independent of $h$.

Finally, we impose {\bf weak consistency}
with the weak problem \eqref{weak_form(a)},
\eqref{weak_form(b)}.  For every pair of sequences
$v^h$, $\psi^h\in H^{1,h}_0$ such that $T^h(v^h)\rightharpoonup v$
weakly in $H^1_0(\Omega)$ and $T^h(\psi^h)\to\psi$ strongly
in $H^1_0(\Omega)$, we have that
\begin{subequations}
\begin{equation}
\label{wk_con_b}
B^h[v^h,\psi^h]\to B[v,\psi]
\end{equation}
and
\begin{equation}
\label{wk_con_f}
F^h(\psi^h)\to F(\psi),
\end{equation}
\end{subequations}
where $B$ and $F$ are defined by \eqref{bilinearform} and
\eqref{weak_form(b)}.

Under these general assumptions, we have the following:

\begin{lem}
\label{abs_convergence}
For every $h>0$, there exists a unique solution $v^h\in H_0^{1,h}$
of the discrete problem
\begin{equation}
\label{discreteproblem}
-B^h[v^h,\psi^h]=F^h(\psi^h),
\end{equation}
for every $\psi^h\in H_0^{1,h}$.

The sequence of extensions $T^h(v^h)$ of the 
family of discrete solutions converges weakly to
$v=u-g+a\chi_{\Omega^-}$ in $H^1_0(\Omega)$, where $u$ the
solution of the weak problem \eqref{weak_form(a)}, 
\eqref{weak_form(b)}.
\end{lem}

\begin{proof}
The existence of a unique solution of the discrete problem
\eqref{discreteproblem} follows by the Lax-Milgram lemma.

By the estimates \eqref{abs_structure(a)}, \eqref{abs_structure(b)},
we obtain a uniform bound for the
sequence of discrete solutions $v^h$
\[
\|v^h\|_h\le C,
\]
with $C$ independent of $h$.  Thus, using \eqref{bdd_ext}
we have the bound 
\[
\|T^h(v^h)\|_{H^1_0(\Omega)}\le C.
\]
By weak compactness in the Hilbert space $H^1_0(\Omega)$,
there is a subsequence $T^{h'}(v^{h'})$ converging weakly to
some $v\in H^1_0(\Omega)$.

Now let $\psi\in H^1_0(\Omega)$ be given.  Using the approximation
property \eqref{str_app}, choose a sequence 
$\psi^{h}\in H^{1,h}_0$ such that
\[
T^{h}(\psi^h)\to \psi
\quad\mbox{in}\quad H^1_0(\Omega).
\]

By weak consistency \eqref{wk_con_b}, \eqref{wk_con_f} we find that the limit
function $v$ satisfies \eqref{weak_form(a)}, \eqref{weak_form(b)}.
Finally, since this problem has a unique solution, it follows
that the full sequence $T^h(v^h)$ converges weakly to $v$
in $H^1_0(\Omega)$.
\end{proof}

\begin{remark}
The simplest way in which to obtain
a consistent, and therefore convergent,
scheme is to define $B^h[v^h,\psi^h]=
B[T^h(v^h),T^h(\psi^h)]$ and $F^h(\psi^h)=F(T^h(\psi^h))$.
This is essentially the method of finite elements
which then hinges on the choice of the extension operator $T^h$,
\cite{ciarlet}.

However, the scheme under consideration in the next sections is
not of this type, insofar as it originates from different approximations
for $B$ and $F$.
\end{remark}

\section{Numerical Method}\label{numerical}

In this section, we rederive the finite difference
scheme from \cite{lfk} for \eqref{pde(a)}, \eqref{pde(b)}, \eqref{pde(c)},
by discretizing the weak formulation \eqref{weak_form(a)},
\eqref{weak_form(b)}.

For the remainder of the paper we assume that the coefficient
matrix is of the form $\beta\;Id$, for some scalar function.
We also assume that the data functions
$a$, $b$, $g$, $f$,
all lie in $C^1(\overline\Omega)$, with $a$, $b$
vanishing on $\partial\Omega$.
For simplicity, 
we restrict ourselves to the special case of a
rectangular domain $\Omega = (x_W,x_E)\times(y_S,y_N)$ in the plane.
Given positive integers $I$ and $J$, set $\Delta x=(x_E-x_W)/(I+1)$
and $\Delta y=(y_N-y_S)/(J+1)$, and
define a uniform grid
$\Omega^h=\{(x_i,y_j)\}$ where
$x_i=x_W+i\Delta x$ and $y_j=y_S+j\Delta y$ for
$i=0,1,\cdots,I+1$ and $j=0,1,\cdots,J+1$.
The grid size is defined as $h=\min(\Delta x,\Delta y)$.  The ratio
$\Delta x/\Delta y$ is fixed when the grid size $h$
goes to zero.

The set of grid functions will be denoted by
\begin{subequations}
\begin{equation}
\label{grid_func}
H^{1,h}=\{w^h=(w_{i,j}):0\le i\le I+1,\; 0\le j\le J+1\}.
\end{equation}
The discrete solution space is defined as
\begin{equation}
\label{discrete_space}
H^{1,h}_0=\{\psi^h=(\psi_{i,j})\in H^{1,h}: 
\psi_{i,j}=0\;\; \mbox{on the grid boundary}\}.
\end{equation}
\end{subequations}

To construct the bilinear form on $H^{1,h}_0$, we discretize
the coefficient $\beta$ in two ways as follows,
\begin{equation}
\label{beta_disc}
\beta^1_{i+1/2,j}=\beta(x_{i+1/2},y_j),
\8\beta^2_{i,j+1/2}=\beta(x_i,y_{j+1/2}).
\end{equation}
For $w^h\in H^{1,h}$
define the usual finite difference operators
\begin{align*}
%\begin{equation}
%\label{fini_dif_op}
%\begin{array}{l}
&(\T_xw)_{i+1/2,j}=(w_{i+1,j}-w_{i,j})/\Delta x,\\
&(\T_yw)_{i,j+1/2}=(w_{i,j+1}-w_{i,j})/\Delta y.
%\end{array}
%\end{equation}
\end{align*}
For $v^h$, $\psi^h\in H^{1,h}_0$, the bilinear form then is given by
\begin{align}
\label{discrete_bilinear}
B^h[v^h,\psi^h]=&
\sum\limits^J_{j=1}\sum\limits^I_{i=0}\beta^1_{i+1/2,j}
(\T_xv)_{i+1/2,j}(\T_x\psi)_{i+1/2,j}\Delta x\Delta y\\
\nonumber
&+\sum\limits^I_{i=1}\sum\limits^J_{j=0}\beta^2_{i,j+1/2}
(\T_yw)_{i,j+1/2}(\T_y\psi)_{i,j+1/2}\Delta x\Delta y.
\end{align}
As in the continuous case, we use this to define a norm
on $H^{1,h}_0$:
\begin{equation}
\label{norm_def}
\|\psi^h\|_h^2=B^h[\psi^h,\psi^h].
\end{equation}

Our next task will be to discretize the linear
functional $F$ in \eqref{weak_form(b)}.
The data functions
naturally give rise to grid functions $a^h$, $b^h$, etc.,
by restriction to the grid $\Omega^h$.
The normal vector is discretized by
\begin{equation}\label{disc_given_func(e)}
n^h=(n^{(1),h},n^{(2),h})=(\phi_x^h,\phi_y^h)/|\nabla \phi^h|
\end{equation}
where for $i=1,\cdots,I$, $j=1,\cdots,J$,
\begin{align*}
  &(\phi_x)_{i,j}=(\phi_{i+1,j}-\phi_{i-1,j})/(2\Delta x),\\
  &(\phi_y)_{i,j}=(\phi_{i,j+1}-\phi_{i,j-1})/(2\Delta y).
\end{align*}
Hence $n^h$ is defined on all interior grid points.
It will not be used on $\partial\Omega$.

We can immediately define 
\begin{subequations}
\begin{equation}
\label{f_1}
F^h_1(\psi^h)=
\sum\limits^I_{i=1}\sum\limits^J_{j=1}f_{i,j}\psi_{i,j}\Delta x\Delta y,
\end{equation}
and
\begin{equation}
\label{f_4}
F^h_2(\psi^h)=B^h[g^h,\psi^h].
\end{equation}

The other two pieces require the localization of integrals
to the subdomain $\Omega^-$.
First, we discretize the characteristic function:
\begin{equation}
\label{chdis}
\chi_{i,j}=\begin{cases}1,&\mbox{if}\quad \phi_{i,j}\le0\\
% (x_{i},y_{j})\in \{\phi\le0\}\\
0,&\mbox{if}\quad \phi_{i,j}>0.
%(x_{i},y_{j})\in \{\phi>0\}.
\end{cases}
\end{equation}
%As was introduced in \cite{lfk},
%$dA_{i+1/2,j}$ will represent
%a consistent approximation 
%of the area of intersection between
%$\Omega^-$ and the two triangular cells sharing the common
%horizontal arm $\{(x,y_j)|x_i\leq x\leq x_{i+1}\}$, shown
%as the shaded region in Figure \ref{fig:2}.
%\begin{figure}
%\centerline{
%\begin{picture}(108,108)
%\put (0,0){\psfig{figure=FIGURES/triangle.eps,height=1.5in,width=1.5in}}
%\put(70,45){\tiny $\theta\Delta x$}
%\put(70,90){$\Omega^-$}
%\put(25,100){$\Omega^+$}
%\put(45,115){$\Gamma$}
%\put(-35,55){\tiny $(x_i,y_j)$}
%\put(-35,100){\tiny $(x_i,y_{j+1})$}
%\put(115,55){\tiny $(x_{i+1},y_j)$}
%\put(115,0){\tiny $(x_{i+1},y_{j-1})$}
%\end{picture}}
%\caption{}
%\label{fig:2}
%\end{figure}
%The area is
%\[dA_{i+1/2,j}=\chi^1_{i+1/2,j}\Delta x\Delta y
%\]
Define  
\begin{align}
\label{chgax}
&\chi^1_{i+1/2,j}
=\left(\chi_{i,j}(1-\theta_{i+1/2,j})+\chi_{i+1,j}\theta_{i+1/2,j}\right)\\
\nonumber
&\theta_{i+1/2,j}=\begin{cases}
|\phi_{i+1,j}|/(|\phi_{i,j}|+|\phi_{i+1,j}|),&
\mbox{if $|\phi_{i,j}|+|\phi_{i+1,j}|>0$,}\\
0&\mbox{otherwise.}
\end{cases}
\end{align}
The factor $\chi^1_{i+1/2,j}\Delta x$ approximates the
length of the portion of the arm
from $(x_i,y_j)$ to $(x_{i+1},y_j)$ within $\Omega^-$.
And also, let
\begin{align}
\label{chgay}
&\chi^2_{i,j+1/2}
=\left(\chi_{i,j}(1-\theta_{i,j+1/2})+\chi_{i,j+1}\theta_{i,j+1/2}\right)\\
\nonumber
&\theta_{i,j+1/2}=\begin{cases}
|\phi_{i,j+1}|/(|\phi_{i,j}|+|\phi_{i,j+1}|),&
\mbox{if $|\phi_{i,j}|+|\phi_{i,j+1}|>0$,}\\
0&\mbox{otherwise.}
\end{cases}
\end{align}

Now we define the remaining two pieces of the linear functional $F$.
\begin{align}
\label{f_2}
F_4^h(\psi^h)=&
\sum\limits^J_{j=1}\sum\limits^I_{i=0}
\T_x(bn^{(1)}\psi)_{i+1/2,j} \,\chi^1_{i+1/2,j}\Delta x\Delta y\\
\nonumber
&+\sum\limits^I_{i=1}\sum\limits^J_{j=0}
\T_y(bn^{(2)}\psi)_{i,j+1/2} \,\chi^2_{i,j+1/2}\Delta x\Delta y
\end{align}
and
\begin{align}
\label{f_3}
F_3^h(\psi^h)=&
-\sum\limits^J_{j=1}\sum\limits^I_{i=0}
\beta^1_{i+1/2,j}(\T_xa)_{i+1/2,j}(\T_x\psi)_{i+1/2,j}\,
\chi^1_{i+1/2,j}\Delta x\Delta y\\
\nonumber
&-\sum\limits^I_{i=1}\sum\limits^J_{j=0}
\beta^2_{i,j+1/2}(\T_ya)_{i,j+1/2}(\T_y\psi)_{i,j+1/2}\,
\chi^2_{i,j+1/2}\Delta x\Delta y.
\end{align}
\end{subequations}

Using \eqref{f_1},$\ldots$,\eqref{f_3}, define
\begin{equation}
\label{dfunc_def}
F^h=F^h_1+\ldots+F^h_4,
\end{equation}
With definitions
\eqref{discrete_bilinear}
and \eqref{dfunc_def}
the discrete problem is  formulated as in \eqref{discreteproblem}.

Next we show that the discrete weak formulation is the
same as the finite difference scheme introduced in \cite{lfk}.

Make the substitution $v^h=u^h-g^h+a^h\chi^h$, using \eqref{chdis}, in 
\eqref{discreteproblem} to write
\begin{multline}
\label{discretesub}
-B^h[u^h,\psi^h]+B^h[g^h,\psi^h]-B^h[a^h\chi^h,\psi^h]\\
=F_1^h(\psi^h)+\cdots+F^h_4(\psi^h).
\end{multline}
We note right away that the second term on the left
cancels with $F^h_2(\psi^h)$ on the right.
We are going to apply summation by parts to remove the difference
operators from the test vector $\psi^h$.  The idea is expressed
by the one-dimensional formula
\begin{equation}
\label{sbp}
-\sum\limits_{i=0}^I\alpha_{i+1/2}(\nabla_x\psi)_{i+1/2}
=\sum\limits_{i=1}^I(\nabla_x\alpha)_i\psi_i,
\end{equation}
provided that $\psi_0=\psi_{I+1}=0$.
Here, $(\nabla_x\alpha)_i=(\alpha_{i+1/2}-\alpha_{i-1/2})/\Delta x$.

Using \eqref{sbp}, the first term in \eqref{discretesub}
can be rewritten as
\begin{subequations}
\begin{equation}
\label{blsbp}
-B^h[u^h,\psi^h]
=\sum\limits_{i=1}^I\sum\limits^J_{j=1}
[\nabla_x(\beta^1\nabla_xu)_{i,j}
+\nabla_y(\beta^2\nabla_yu)_{i,j}]\psi_{i,j}
\Delta x\Delta y.
\end{equation}

In the same way, we have using definitions \eqref{chgax},
\eqref{chgay},
\begin{equation}
\label{f2sbp}
F^h_4(\psi^h)=
-\sum\limits_{i=1}^I\sum\limits^J_{j=1}
b_{i,j}[n^{(1)}_{i,j}(\nabla_x\chi^1)_{i,j}
+n^{(2)}_{i,j}(\nabla_y\chi^2)_{i,j}]\psi_{i,j}
\Delta x\Delta y
\end{equation}

To treat the remaining term on the left-hand side
of \eqref{discretesub}, we use the following product
rule for the difference operator
\begin{align*}
&\nabla_x(a\chi)_{i+1/2,j}=(\nabla_xa)_{i+1/2,j}\chi^1_{i+1/2,j}
+a^1_{i+1/2,j}(\nabla_x\chi)_{i+1/2,j}\\
&a^1_{i+1/2,j}=a_{i+1,j}(1-\theta_{i+1/2,j})+a_{i,j}\theta_{i+1/2,j},
\end{align*}
in which $\chi^1_{i+1/2,j}$ and $\theta_{i+1/2,j}$
were defined in \eqref{chgax}.
Similarly, we have using \eqref{chgay}
\begin{align*}
&\nabla_y(a\chi)_{i,j+1/2}=\nabla_ya_{i,j+1/2}\chi^2_{i,j+1/2}
+a^2_{i,j+1/2}\nabla_x\chi_{i,j+1/2}\\
&a^2_{i,j+1/2}=a_{i,j+1}(1-\theta_{i,j+1/2})+a_{i,j}\theta_{i,j+1/2}.
\end{align*}
It follows from this and \eqref{sbp} that
\begin{align}
\nonumber
B^h[a^h\chi^h,\psi^h]+&F^h_3(\psi^h)\\
\nonumber
=&\sum\limits_{i=0}^I\sum\limits^J_{j=1}
\beta^1_{i+1/2,j}a^1_{i+1/2,j}
(\T_x\chi)_{i+1/2,j}(\T_x\psi)_{i+1/2,j}
\Delta x\Delta y\\
\label{f3ibp}
&+\sum\limits_{i=1}^I\sum\limits^J_{j=0}
\beta^2_{i,j+1/2}a^2_{i,j+1/2}
(\T_y\chi)_{i,j+1/2}(\T_y\psi)_{i,j+1/2}
\Delta x\Delta y\\
=&-\sum\limits_{i=1}^I\sum\limits^J_{j=1}
[\T_x(\beta^1 a^1\T_x\chi)_{i,j}
+\T_y(\beta^2 a^2\T_x\chi)_{i,j}]\psi_{i,j}
\Delta x\Delta y
\nonumber
\end{align}
\end{subequations}

Combining \eqref{discretesub}, \eqref{blsbp}, \eqref{f2sbp},
\eqref{f3ibp},
we obtain
\begin{align*}
\sum\limits_{i=1}^I\sum\limits^J_{j=1}&
[\nabla_x(\beta^1\nabla_xu)_{i,j}
+\nabla_y(\beta^2\nabla_yu)_{i,j}]\psi_{i,j}
\Delta x\Delta y\\
=&-\sum\limits_{i=1}^I\sum\limits^J_{j=1}
[\T_x(\beta^1 a^1\T_x\chi)_{i,j}
+\T_y(\beta^2 a^2\T_x\chi)_{i,j}]\psi_{i,j}
\Delta x\Delta y\\
&-\sum\limits_{i=1}^I\sum\limits^J_{j=1}
b_{i,j}[n^{(1)}_{i,j}\nabla_x(\chi^1)_{i,j}
+n^{(2)}_{i,j}\nabla_y(\chi^2)_{i,j}]\psi_{i,j}
\Delta x\Delta y\\
&+\sum\limits_{i=1}^I\sum\limits^J_{j=1}
f_{i,j}\psi_{i,j}
\Delta x\Delta y.
\end{align*}
Since this must hold for all test vectors $\psi^h\in H^{1,h}_0$,
we have shown that the finite difference scheme
\begin{multline*}
\nabla_x(\beta^1\nabla_xu)_{i,j}
+\nabla_y(\beta^2\nabla_yu)_{i,j}
=-\T_x(\beta^1 a^1\T_x\chi)_{i,j}
-\T_y(\beta^2 a^2\T_x\chi)_{i,j}\\
-b_{i,j}[n^{(1)}_{i,j}(\nabla_x\chi^1)_{i,j}
+n^{(2)}_{i,j}(\nabla_y\chi^2)_{i,j}]
+f_{i,j}
\end{multline*}
holds at all interior grid points.
Note that this is the scheme that was 
found in \cite{lfk}, see equation (77) therein.

\begin{remark}
Here the passage from the weak formulation via
summation by parts to the
finite difference scheme is analogous to what is
often done with PDE's.
%Conversely, by multiplying a finite difference scheme by
%a discrete test 
%function, summing, and using summation by parts, 
%we get discrete weak problems from finite difference methods.
\end{remark}

%\begin{remark}
%New finite element methods can be directly derived from
%the weak formulation (\ref{weak_form(a)},\ref{weak_form(b)}). 
%Because the interface conditions are on the right 
%handside of the weak formulation 
%(\ref{weak_form(a)},\ref{weak_form(b)}), the basis functions can be continuous
%as one usually does for problems without interfaces.
%\end{remark}

\section{Convergence}

In this section, we establish the converge of the
scheme \eqref{discreteproblem}.
%with the bilinear form \eqref{discrete_bilinear}
%and linear functional \eqref{dfunc_def}
%defined in the previous section.  The strategy
%will be to apply Lemma \ref{abs_convergence}.
%In order to do so, we must verify the hypotheses:
%uniform boundednes of $B^h[\cdot,\cdot]$ and $F^h(\cdot)$,
%the existence of a family of extension operators
%with the strong approximation property, and
%finally the weak consistency of the problem.

\begin{thm}
\label{main_result}
Let $\Omega\subset\rr^2$ be a rectangle. 
Assume that the data functions $a$, $b$, $g$, $f$
all lie in $C^1(\overline\Omega)$, with $a$,
$b$ vanishing on $\partial\Omega$.
Suppose that the coefficients have the form $\beta\; Id$.
Then there is a family of linear
extensions $T^h:H^{1,h}_0\to H^1_0(\Omega)$ which
together with the bilinear forms $B^h$ \eqref{discrete_bilinear},
and the linear functionals $F^h$ \eqref{dfunc_def}
satisfy the structural conditions
\eqref{bdd_ext}, \eqref{str_app}, \eqref{abs_structure(a)},
\eqref{abs_structure(b)}, \eqref{wk_con_b}, \eqref{wk_con_f}.

The sequence of 
extended approximate solutions $T^h(v^h)$ of the discrete weak problem
\eqref{discreteproblem}, \eqref{discrete_bilinear}, \eqref{dfunc_def}
converge weakly to the weak solution of the PDE \eqref{weak_form(a)},
\eqref{weak_form(b)} in $H^1_0(\Omega)$.
\end{thm}

\begin{proof}
The second statement follows from the first by Lemma
\ref{abs_convergence}.
The next three subsections are devoted to the verification
of the structural conditions:  uniform boundedness,
extension and approximation, and finally consistency.
\end{proof}

\subsection{Uniform boundedness}
\begin{lem}
\label{uniform_bilinear}
The family of bilinear forms \eqref{discrete_bilinear}
satisfy the uniform bounds \eqref{abs_structure(a)}.
\end{lem}
\begin{proof}
Because of our choice of norm in \eqref{norm_def},
the lower bound in \eqref{abs_structure(a)} is immediate.
The upper bound in \eqref{abs_structure(a)} follows
easily from the Cauchy-Schwarz inequality.
\end{proof}

\begin{lem}
\label{uniform_func}
The family of linear functionals \eqref{dfunc_def}
satisfy the uniform bound \eqref{abs_structure(b)}.
\end{lem}
\begin{proof}
We treat the four pieces $F^h_1,\ldots,F^h_4$ in turn.

By the Cauchy-Schwarz inequality, we have
\[
|F^h_1(\psi^h)|\le (\sum_{i,j}|f_{i,j}|^2 \Delta x\Delta y)^{1/2}
(\sum_{i,j}|\psi_{i,j}|^2 \Delta x\Delta y)^{1/2}.
\]
The first factor is just a Riemann sum for the $L^2$-norm
of $f$, and so it is bounded by $2\|f\|_{L^2(\Omega)}$,
for $h$ small enough.  The second factor is estimated
using the discrete version of the Poincar\'e inequality
\begin{equation}
\label{poincare}
(\sum_{i,j}|\psi_{i,j}|^2 \Delta x\Delta y)^{1/2}
\le C(\Omega)(\sum\limits_{i=0}^I\sum\limits_{j=1}^J
|(\nabla_x\psi)_{i+1/2,j}|^2 \Delta x\Delta y)^{1/2},
\end{equation}
which follows as in the continuous case using summation by parts.
This last sum is then bounded by $m^{-1}\|\psi^h\|_h$ because
of the uniform lower bound for the coefficients $\beta$.

From \eqref{f_2},
we can use the discrete product formula to rewrite $F^h_4$ as
\begin{multline*}
F^h_4(\psi^h)\\
= \sum\limits_{i=0}^I\sum\limits_{j=1}^J
[(bn^{(1)})_{i+1,j}(\nabla_x\psi)_{i+1/2,j}
+\nabla_x(bn^{(1)})_{i+1/2,j}\psi_{i,j}]\chi^1_{i+1/2,j}\Delta x\Delta y\\
+\sum\limits_{i=1}^I\sum\limits_{j=0}^J
[(bn^{(2)})_{i,j+1}(\nabla_x\psi)_{i,j+1/2}
+\nabla_x(bn^{(2)})_{i,j+1/2}\psi_{i,j}]\chi^2_{i,j+1/2}\Delta x\Delta y
\end{multline*}
Since $b$ and $(n^{(1)},n^{(2)})$ are in $C^1$ and $(\chi^1,\chi^2)$ is
bounded,
using the Cauchy-Schwarz inequality again, we get the bound
\[
|F^h_4(\psi^h)|\le C\|bn\|_{L^2(\Omega)} \|\psi^h\|_h
+C\|\nabla\cdot(bn)\|_{L^2(\Omega)}
\left(\sum_{i,j} |\psi_{i,j}|^2\Delta x\Delta y\right)^{1/2}.
\]
Applying \eqref{poincare} to the last sum above,
we obtain the desired bound for $F^h_4$.

In the same way, we get,
\[
|F^h_2(\psi^h)|\le C \|\beta\nabla g\|_{L^2(\Omega)}\|\psi^h\|_h
\5\mbox{and}\5|F^h_3(\psi^h)|\le C \|\beta\nabla a\|_{L^2(\Omega)}\|\psi^h\|_h.
\]
%and
%\[
%|F^h_3(\psi^h)|\le C \|\beta\nabla a\|_{L^2(\Omega)}\|\psi^h\|_h.
%\]
\end{proof}

\subsection{Strong Approximation}

In this section we define a uniformly
bounded family of extension operators
$T^h$ using the 
basic approach from the theory of finite elements.
Then we verify the strong approximation property.

\begin{figure}\label{fig:1}
\begin{picture}(108,108)
\put(120,0){\psfig{figure=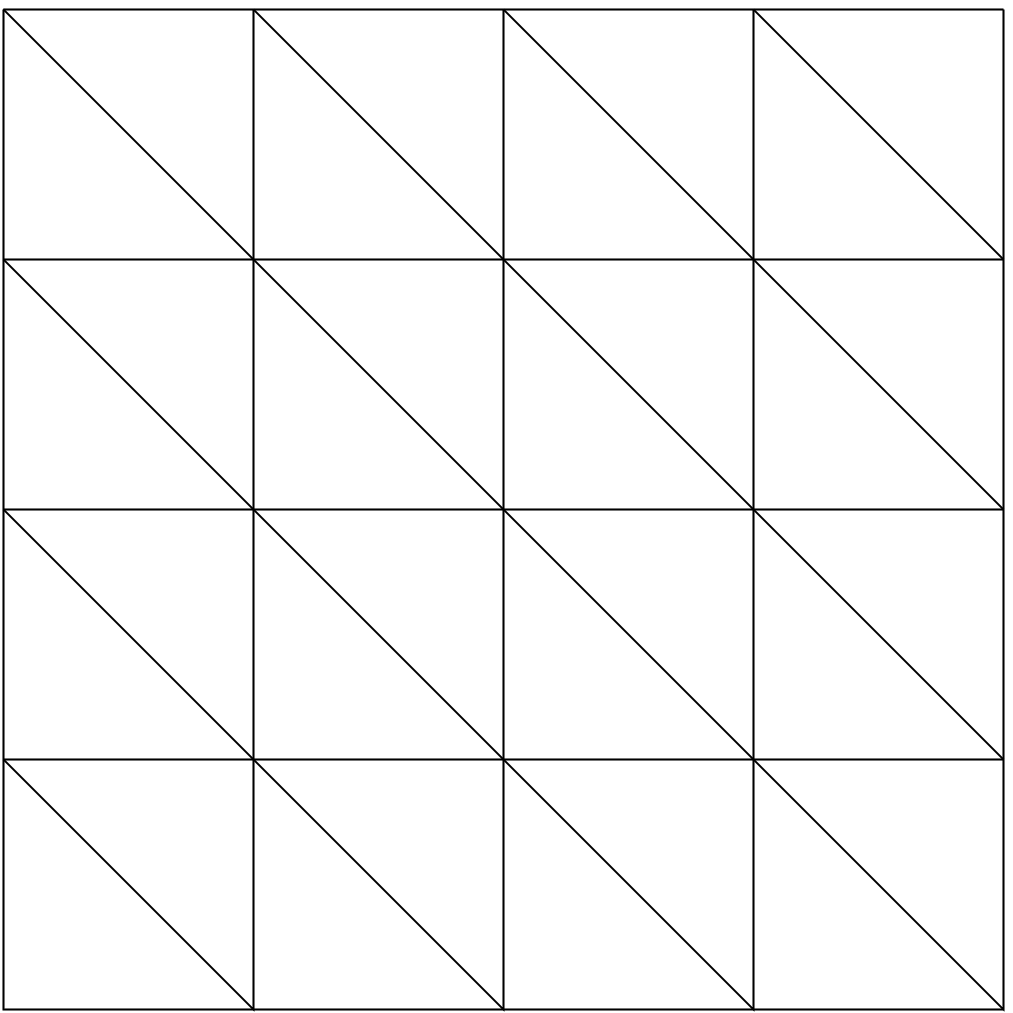,height=1.5in,width=1.5in}}
\end{picture}
\caption{}
\end{figure}

To this end, for each $h$ consider a triangulation of
the domain $\Omega$ containing all triangles with vertices
\[
\{(x_i,y_j),(x_{i+1},y_j),(x_i,y_{j+1})\}
\quad \mbox{or}\quad
\{(x_i,y_j),(x_{i-1},y_j),(x_i,y_{j-1})\}
\]
based on the grid $\Omega^h$, as
shown in Figure \ref{fig:1}.  For any 
grid point $(x_i,y_j)\in \Omega^h$,
let $\eta_{i,j}^h\in H^1(\Omega)$
be the continuous, piecewise linear function which
is equal to 1 at the grid point $(x_i,y_j)$ and
equal to 0 at all other grid points.
Given $w^h\in H^{1,h}$,
the extension operator is defined as
\[
T^h(w^h)=\sum_{i,j}w_{i,j}\eta_{i,j}.
\]
Then $T^h:H^{1,h}\to H^1(\Omega)$ and $T^h:H^{1,h}_0\to H^1_0(\Omega)$.
Although it is not required by the abstract framework,
the extensions $T^h$ are linear operators.

Since $T^h(w^h)$ is the unique continuous linear interpolant
on each triangle $\{(x_i,y_j),(x_{i\pm1},y_j),(x_i,y_{j\pm1})\}$,
we have explicitly
\begin{equation}
\label{interpolant}
T^h(w^h)(x,y)=w_{i,j}+(\T_xw)_{i\pm1/2,j}(x-x_{i\pm1})
+(\T_yw)_{i,j\pm1/2}(y-y_{j\pm1}).
\end{equation}

Given the coefficient function $\beta$ on $\Omega$, we
write using the definition \eqref{beta_disc}
\[
\beta^{1,h}=(\beta^1_{i+1/2,j}),\quad
\beta^{2,h}=(\beta^2_{i,{j+1/2}}).
\]

In general, given a discrete function $\beta^{1,h}$ defined at the
half grid points $(x_{i+1/2},y_j)$, as above,
we define the piecewise constant extension
\[
S^{1,h}(\beta^{1,h})(x,y)=\beta^1_{i+1/2,j}
\]
on every triangular cell having 
the horizontal edge from $(x_i,y_j)$ to $ (x_{i+1},y_j)$, see Figure 2.
Similarly, given a discrete function $\beta^{2,h}$ defined at the
half grid points $(x_i,y_{j+1/2})$,
we define
\[
S^{2,h}(\beta^{2,h})(x,y)=\beta^2_{i,j+1/2}
\]
on every triangular cell having the vertical edge
from $(x_i,y_j)$ to $ (x_i,y_{j+1})$, see Figure 2.

\begin{figure}\label{fig:2}
\begin{picture}(200,200)
\put(54,0){\psfig{figure=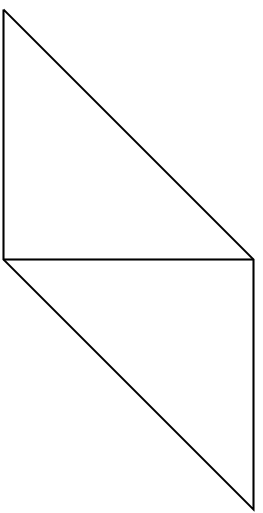,height=2.4in,width=1.2in}}
\put(25,85){\tiny $(x_i,y_j)$}
\put(145,85){\tiny $(x_{i+1},y_j)$}
\put(78,96){\tiny $(x_{i+1/2},y_j)$}
\put(96,84){\small $\bullet$}
\put(208,25){\psfig{figure=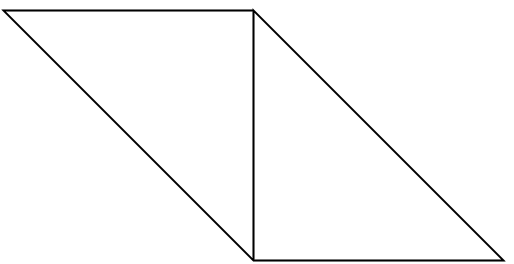,height=1.2in,width=2.4in}}
\put(292,70){\small $\bullet$}
\put(250,78){\tiny $(x_i,y_{j+1/2})$}
\put(282,17){\tiny $(x_i,y_j)$}
\put(277,119){\tiny $(x_i,y_{j+1})$}
\end{picture}
\caption{}
\end{figure}

The key to our estimates will be to
replace sums by integrals.  The following
lemma summarizes the important formulas.

\begin{lem}
\label{disc_to_cont}
With the definitions above, we have
\begin{align}
\label{d2c_bil}
&B^h[v^h,\psi^h]=
\int_\Omega\Bigl[S^{1,h}(\beta^{1,h})T^{h}(v^h)_xT^h(\psi^h)_x\\
&\hspace{1in}
\nonumber
+S^{2,h}(\beta^{2,h})T^h(v^h)_yT^h(\psi^h)_y\Bigr]
\\
\label{d2c_f2}
&F^h_4(\psi^h)=\int_\Omega\Bigl[S^{1,h}(\chi^{1,h})
T^h(b^hn^{(1),h}\psi^h)_x\\
&\hspace{1in}
\nonumber
+S^{2,h}(\chi^{2,h})T^h(b^hn^{(2),h}\psi^h)_y\Bigr]
\\
\label{d2c_f3}
&F^h_3(\psi^h)=
-\int_\Omega\Bigl[S^{1,h}(\beta^{1,h}) S^{1,h}(\chi^{1,h})
T^{h}(a^h)_xT^h(\psi^h)_x \\
&\hspace{1in}
\nonumber
+S^{2,h}(\beta^{2,h}) S^{2,h}(\chi^{2,h})
T^h(a^h)_yT^h(\psi^h)_y\Bigr].
\end{align}
%Let $\psi^h\in H^{1,h}_0$.  Suppose $\beta^{1,h}$, $\beta^{2,h}$
%are defined at the half grid points.
%Then 
%\[
%\sum\limits_{i=0}^I\sum_{j=1}^J
%\beta^1_{i+1/2,j}(\T_x\psi)_{i+1/2,j}\Delta x\Delta y
%=\int_\Omega S^{1,h}(\beta^{1,h})T^h(\psi^h)_x,
%\]
%and
%\[
%\sum\limits^I_{i=1}\sum\limits^J_{j=0}
%\beta^2_{i,j+1/2}(\T_y\psi)_{i,j+1/2}\Delta x\Delta y=
%  \int\limits_\Omega S^{1,h}(\beta^{2,h})T^h(\psi^h)_y.
%\]
\end{lem}

\begin{proof}
Suppose that $\xi^1_{i+1/2,j}$ is defined for $i=0,\ldots,I$,
$j=0,\ldots,J+1$, and $\xi^1_{i+1/2,j}=0$ for $j=0$ or $j=J+1$.
Each point of the form $(x_{i+1/2},y_j)$,
$i=1,\ldots,I$, $j=1,\ldots,J$, lies on the boundary of
the two triangular cells sharing the arm from
$(x_i,y_j)$ to $(x_{i+1},y_j)$ , and on these cells
we have 
$\xi^1_{i+1/2,j}=S^{1,h}(\xi^{1,h})$(x,y).
On the remaining cells, namely those with $j=0$ or $j=J+1$,
we have $\xi^1_{i+1/2,j}=S^{1,h}(\xi^{1,h})=0$.
Therefore, since the area of each pair of cells (where
$S^{1,h}(\xi^{1,h})$ could be nonzero)
is equal to $\Delta x\Delta y$, we have
\[
\sum\limits_{i=0}^I\sum_{j=1}^J \xi^1_{i+1/2,j}\Delta x\Delta y
=\int_\Omega S^{1,h}(\xi^{1,h}).
\]
Similarly, we have
\[
\sum\limits_{i=1}^I\sum_{j=0}^J \xi^2_{i,j+1/2}\Delta x\Delta y
=\int_\Omega S^{2,h}(\xi^{1,h}),
\]
provided that $\xi^2_{i,j+1/2}=0$ for $i=0$ or $i=I+1$.

Notice that the piecewise continuous extensions $S^{\alpha,h}$
are multiplicative in the sense that, for $\alpha=1,2$, 
$S^{\alpha,h}(a^{\alpha,h}b^{\alpha,h})
=S^{\alpha,h}(a^{\alpha,h})S^{\alpha,h}(b^{\alpha,h})$,
for arbitrary $a^{\alpha,h}$, $b^{\alpha,h}$ defined
at the half grid points.
So using the formulas just derived, we have from
\eqref{discrete_bilinear} that
\begin{align*}
B^h[v^h,\psi^h]
&= \int_\Omega[S^{1,h}(\beta^{1,h}\T_xv^h\T_x\psi^h)
+S^{2,h}(\beta^{2,h}\T_yv^h\T_y\psi^h)]\\
&= \int_\Omega[S^{1,h}(\beta^{1,h})S^{1,h}(\T_xv^h)S^{1,h}(\T_x\psi^h)\\
&\mbox{\hspace{1in}}
+S^{2,h}(\beta^{2,h})S^{2,h}(\T_yv^h)S^{2,h}(\T_y\psi^h)].
\end{align*}
By \eqref{interpolant}, we see that
$S^{1,h}(\T_x\psi^h)=T^h(\psi^h)_x$
and $S^{2,h}(\T_y\psi^h)=T^h(\psi^h)_y$,
and therefore, we have verified \eqref{d2c_bil}.

The proofs of the other two formulas are similar
and will be omitted.
\end{proof}

We are now ready to study the extensions.

\begin{lem}
\label{ext_lem_1}
The family of extensions $T^h:H^{1,h}_0\to H^1_0(\Omega)$ 
satisfy the uniform estimates
\[
\frac{m}{M}
\|\psi^h\|_h^2
\le
\|T^h\psi^h\|_{H^1_0(\Omega)}^2
\le \frac{M}{m}
\|\psi^h\|_h^2
\]
\end{lem}

\begin{proof}
Let $\psi^h\in H^{1,h}_0$.
By definition, we have
\[
\|T^h\psi^h\|_{H^1_0(\Omega)}^2
=B[T^h\psi^h,T^h\psi^h]
\quad\mbox{and}\quad
\|\psi^h\|_h^2=B^h[\psi^h,\psi^h].
\]
Moreover, by \eqref{d2c_bil}, we have that
\[
B^h[\psi^h,\psi^h]=
\int_\Omega[S^{1,h}(\beta^{1,h})T^h(\psi^h)_x^2
+S^{2,h}(\beta^{2,h})T^h(\psi^h)_y^2].
\]

Now by the uniform bounds on $\beta$, we have that
\[
\beta(x,y)\le \frac{M}{m} \beta(x',y'),
\]
for arbitrary $(x,y)$, $(x',y')\in\Omega$,
and so
\begin{equation}
\label{uplo}
\frac{m}{M}\beta\le
S^{1,h}(\beta^{1,h}), S^{2,h}(\beta^{2,h})\le \frac{M}{m}\beta.
\end{equation}

By \eqref{uplo}, it follows that
\[
\frac{m}{M} B^h[\psi^h,\psi^h]\le
B[T^h\psi^h,T^h\psi^h]
\le \frac{M}{m} B^h[\psi^h,\psi^h].
\]
\end{proof}

\begin{lem}
\label{ext_lem_2}
The extension operator $T^h:H^{1,h}_0\to H^1_0(\Omega)$
has the strong approximation property \eqref{str_app}.
\end{lem}

\begin{proof}
First let $\psi\in C^1_0(\overline{\Omega})$, and define $\psi^h\in H^{1,h}_0$
to be the grid function with values $\psi_{i,j}=\psi(x_i,y_j)$.
Since $\psi\in C^1(\overline{\Omega})$, $\psi$ and its
first derivatives are uniformly continuous on $\overline{\Omega}$.
It follows from
\eqref{interpolant}
the piecewise constant functions $T^h(\psi^h)_x$
and $T^h(\psi^h)_y$ converge uniformly to 
$\psi_x$ and $\psi_y$, respectively, therefore
also in $L^2(\omega)$:
\[
\|\nabla(T^h\psi^h-\psi)\|_{L^2(\Omega)}\to0,
\]
as $h\to0$.
Since $\beta$ has a uniform lower bound,
this implies that
\[
\|T^h\psi^h-\psi\|_{H^1_0(\Omega)}\to0,
\]
as $h\to0$.

Finally, the result for general $\psi\in H^1_0(\Omega)$
follows by density.
\end{proof}

\subsection{Weak consistency}

\begin{lem}
\label{bl_con_lem}
The bilinear form \eqref{discrete_bilinear}
satisfies the weak consistency hypothesis \eqref{wk_con_b}.
\end{lem}

\begin{proof}
Suppose that $v^h$, $\psi^h\in H^{1,h}_0$ are sequences
such that $T^h(v^h)\rightharpoonup v$ weakly in $H^1_0(\Omega)$
and $T^h(\psi^h)\to \psi$ strongly in $H^1_0(\Omega)$.
Thus, we have that $\T T^h(v^h) \rightharpoonup \T v$
weakly in $L^2(\Omega)$ and $\T T^h(\psi^h)\to \T\psi$
strongly in $L^2(\Omega)$.

Recalling the definitions 
from the previous subsection, we have that the functions
$S^{\alpha,h}(\beta^{\alpha,h})$ ($\alpha=1,2$) are both
uniformly bounded and converge pointwise to
$\beta$ in $\Omega\setminus\Gamma$.  
Writing 
\[
S^{1,h}(\beta^{1,h}) T^h(\psi^h)_x-\beta\psi_x
=[S^{1,h}(\beta^{1,h})-\beta]\psi_x
+S^{1,h}(\beta^{1,h})[ T^h(\psi^h)_x-\psi_x],
\]
it follows by the dominated convergence theorem
that $S^{1,h}(\beta^{1,h})T^h(\psi^h)_x\to\beta\psi_x$
strongly in $L^2(\Omega)$.  

Now, using Lemma \eqref{d2c_bil} we have that
\begin{align*}
B^h[v^h,\psi^h]
&=\int_\Omega[S^{1,h}(\beta^{1,h})T^{h}(v^h)_xT^h(\psi^h)_x\\
&\mbox{\hspace{1in}}
+S^{2,h}(\beta^{2,h})T^h(v^h)_yT^h(\psi^h)_y]\\
&=\langle T^h(v^h)_x, 
S^{1,h}(\beta^{1,h}) T^h(\psi^h)_x \rangle_{L^2(\Omega)}\\
&\mbox{\hspace{1in}}
+\langle T^h(v^h)_y, 
S^{2,h}(\beta^{2,h}) T^h(\psi^h)_y \rangle_{L^2(\Omega)}\\
&\to \langle  v_x, \beta\psi_x \rangle_{L^2(\Omega)}
+\langle v_y, \beta \psi_y \rangle_{L^2(\Omega)}\\
&=B[v,\psi],
\end{align*}
as $h\to0$.  
\end{proof}

\begin{lem}
The linear functional $F$ defined in \eqref{dfunc_def}
satisfies the weak consistency hypothesis \eqref{wk_con_f}.
\end{lem}

\begin{proof}
Let $\psi\in H^1_0(\Omega)$ and suppose that $\psi^h\in H^{1,h}_0$
is a sequence such that $T^h(\psi^h)\to \psi$ strongly in
$H^1_0(\Omega)$.
We must show that $F^h(\psi^h)\to F(\psi)$.

To begin, we observe that it is enough to prove this
for test functions in $C^1_0(\overline\Omega)$.
For any $\overline{\psi}\in C^1_0(\overline\Omega)$, we can write
\begin{multline*}
|F^h(\psi^h)-F(\psi)|\le |F^h(\psi^h-\overline{\psi}^h)|
+|F^h(\overline{\psi}^h)-F(\overline\psi)|
+|F(\overline{\psi}-\psi)|\\
\le\|F^h\|\|\psi^h-\overline{\psi}^h\|_h
+|F^h(\overline{\psi}^h)-F(\overline\psi)|
+\|F\|\;\|\overline{\psi}-\psi\|_{H^1_0(\Omega)}.
\end{multline*}
By Lemma \ref{ext_lem_1},
we have
\begin{multline*}
\|\psi^h-{\overline\psi}^h\|_h
\le \frac{M}{m}\|T^h(\psi^h-\overline{\psi}^h)\|_{H^1_0(\Omega)}\\
\le \frac{M}{m}
\left[\|T^h(\psi^h)-\psi\|_{H^1_0(\Omega)}
+\|T^h(\overline{\psi}^h)-\overline{\psi}\|_{H^1_0(\Omega)}
+\|\overline{\psi}-\psi\|_{H^1_0(\Omega)}\right].
\end{multline*}
By Lemma \ref{uniform_func}, 
the norms $\|F^h\|$, $\|F\|$ are uniformly bounded.  
Therefore, we have shown that
\begin{multline*}
|F^h(\psi^h)-F(\psi)|\le
|F^h(\overline{\psi}^h)-F(\overline\psi)|
+C\Big[ \|T^h(\psi^h)-\psi\|_{H^1_0(\Omega)}\\
+\|T^h(\overline{\psi}^h)-\overline{\psi}\|_{H^1_0(\Omega)}
+\|\overline{\psi}-\psi\|_{H^1_0(\Omega)}\Big].
\end{multline*}
We claim that this can be made arbitrarily small for all $h<h_0$.
By assumption,
we have that $T^h(\psi^h)-\psi\to0$ in ${H^1_0(\Omega)}$.
By density, $\overline\psi$ can be chosen arbitrarily
close to $\psi$ in ${H^1_0(\Omega)}$.
By construction, we have that
$T^h(\overline{\psi}^h)-\overline{\psi}\to0$ in ${H^1_0(\Omega)}$.
This covers all but the first term above.
The first term can also be made small if the consistency
condition  is valid for $\overline\psi\in C^1_0(\overline\Omega)$.

We now proceed to verify the consistency of $F$ under
the assumption that $\psi\in C^1_0(\overline\Omega)$
by considering each individual piece.

By \eqref{f_1} and the
fact that $\psi\in C^1_0(\overline\Omega)$, 
we see that $F_1^h(\psi^h)$ is simply
a Riemann sum for $F_1(\psi)$, and thus
$F^h_1(\psi^h)\to F_1(\psi)$, as $h\to0$.

Given $g\in C^1(\overline\Omega)$, we have 
that 
$T^h(g^h)\to g$ in $H^1_0(\Omega)$,
as in the proof of Lemma \ref{ext_lem_2}.
Thus, we may apply the result of Lemma \ref{bl_con_lem}
to conclude that
\[
F_2(\psi^h) = B^h[g^h,\psi^h] \to B[g,\psi] = F_4(\psi).
\]

According to \eqref{d2c_f2}, we have
\begin{multline*}
F^h_4(\psi^h)=\int_\Omega[S^{1,h}(\chi^{1,h})
T^h(b^hn^{(1),h}\psi^h)_x
+S^{2,h}(\chi^{2,h})T^h(b^hn^{(2),h}\psi^h)_y].
\end{multline*}
Since $b$, $n$, $\psi\in C^1(\overline\Omega)$,
we have as in the proof of Lemma \ref{ext_lem_2},
that $\nabla T^h(b^hn^h\psi^h)\to \nabla bn\psi$ in $L^2$.
Moreover, $S^{\alpha,h}(\chi^{\alpha,h})$ is uniformly
bounded and tends to $\chi_{\Omega^-}$ pointwise
in $\Omega\setminus\Gamma$.  Thus, it follows that
$F^h(\psi^h)\to F(\psi)$.

The consistency of $F_3^h$ follows from the formula
\eqref{d2c_f3}, in the same way.
\end{proof}

{\bf Acknowledgment:} We thank S.\ Osher and R.\ Fedkiw for their
stimulating discussions.  In particular S.\ Osher
gave us the insightful suggestion to try the weak formulation.

\end{document}